\documentclass[12pt,twoside]{amsart}
\usepackage{amssymb}
\usepackage{graphicx}
\nonstopmode

\numberwithin{equation}{section}
\newcommand{\NN}{\mathbb{N}}
\newcommand{\PP}{\mathbb{P}}

\newcommand{\ZZ}{\mathbb{Z}}
\newcommand{\HC}{\mathcal{H}}

\newcommand{\al}{\alpha}
\newcommand{\be}{\beta}
\newcommand{\ga}{\gamma}
\newcommand{\de}{\delta}
\newcommand{\la}{\lambda}

\DeclareMathOperator{\pnt}{\raise 0.5mm \hbox{\large\bf.}}

\DeclareMathOperator{\charf}{char}

\def\cocoa{{\hbox{\rm C\kern-.13em o\kern-.07em C\kern-.13em o\kern-.15em A}}}

\newtheorem{theorem}{Theorem}[section]
\newtheorem{lemma}[theorem]{Lemma}
\newtheorem{proposition}[theorem]{Proposition}
\newtheorem{corollary}[theorem]{Corollary}
\newtheorem{conjecture}[theorem]{Conjecture}

\theoremstyle{definition}

\newtheorem{remark}[theorem]{Remark}

\newtheorem{notation}[theorem]{Notation}

\begin{document}

\title[The weak Lefschetz property, monomial ideals, and lozenges]{The weak Lefschetz property, monomial ideals,\\ and lozenges}

\author[D.\ Cook II, U.\ Nagel]{David Cook II, Uwe Nagel}
\address{Department of Mathematics, University of Kentucky,
715 Patterson Office Tower, Lexington, KY 40506-0027, USA}
\email{dcook@ms.uky.edu, uwenagel@ms.uky.edu}
\thanks{Part of the work for this paper was done while the authors were partially supported by the National Security Agency under Grant Number H98230-09-1-0032.}


\begin{abstract}
    We study the weak Lefschetz property and the Hilbert function of level Artinian monomial almost complete intersections in three variables.  Several such
    families are shown to have the weak Lefschetz property if the characteristic of the base field is zero or greater than the maximal degree of any minimal
    generator of the ideal.  Two of the families have an interesting relation to tilings of hexagons by lozenges.  This lends further evidence to a
    conjecture by Migliore, Mir\'o-Roig, and the second author.  Finally, using our results about the weak Lefschetz property, we show that the Hilbert function
    of each level Artinian monomial almost complete intersection in three variables is peaked strictly unimodal.
\end{abstract}

\maketitle


%
%
\section{Introduction} \label{sec:introduction}

Let $A$ be a standard graded Artinian algebra over a field $K$.  Then $A$ has the \emph{weak Lefschetz property} if there is a linear form $\ell \in A$ such
that, for all integers $d$, the multiplication map
$$ \times\ell: [A]_d \rightarrow [A]_{d+1} $$
has maximal rank, that is, it is surjective or injective.  In this case, the linear form $\ell$ is called a \emph{Lefschetz element} of $A$.

This property is of interest mainly because it constrains the Hilbert function as shown in~\cite{HMNW}, which in turn has interesting consequences
(see, e.g., \cite{St} for a spectacular application). Furthermore, it is a difficult task to classify which rings do (and do not) have the weak
Lefschetz property. For example, in~\cite{HMNW} it was shown that all height three complete intersections over a field of characteristic zero have
the weak Lefschetz property, but this is still unknown if we  consider height four complete intersections.

In this note we further explore level Artinian monomial almost complete intersections in three variables, as discussed in~\cite{Br},~\cite{BK}, and
more extensively in~\cite{MMN}.  Even in this restricted setting, it is still unclassified which rings have the weak Lefschetz property. However, in
\cite{MMN} a conjectural solution is put forth, restated here along with known partial results in Section~\ref{sec:conj}.

Several parts of this conjecture  have been established in \cite{MMN}. Here we resolve some of the open cases, thus lending further evidence to the
conjecture.  In Section~\ref{sec:simple-cases} we consider three rather straightforward cases, where three of the four parameters are equal.
In Section~\ref{sec:extremal-cases} we consider the two cases where a parameter is extremal.

The key to these results is the computation of a certain determinant which was shown  to play a crucial role in \cite{MMN}. Interestingly, the
computation of the determinant in the two extremal cases reveals a connection to combinatorial objects, namely to tilings of hexagons by lozenges.

While the conjecture in \cite{MMN} is for algebras over fields of characteristic zero only, our computation of the determinants allows us also to
establish the weak Lefschetz property also over fields of sufficiently large characteristic.  In fact, we give an effective lower bound on the
characteristic in each case.  However, in Remark~\ref{r:extremal-WLP} we notice that in general the maximal degree of the minimal generators gives
no indication of a such bound on the characteristic.

Last, in Section~\ref{sec:unimodal} we show,  using also our results from Section~\ref{sec:extremal-cases}, that \emph{every} level Artinian monomial
almost complete intersection $R/I$ has a \emph{peaked strictly unimodal} Hilbert function; that is, if $h$ is the Hilbert function of $R/I$, then
\[ h(0) < h(1) < \cdots < h(s) = h(s+1) = \cdots = h(s+t-1) > h(s+t) > \cdots > h(e) \]
where $s, \ldots, s+t-1$ are the peak degrees and $e$ is the socle degree of $R/I$.  This result in turn gives a partial answer to Question~8.2(1)
from~\cite{MMN}. It shows that for these algebras the knowledge of the Hilbert function does not allow one to decide whether the algebra has the weak
Lefschetz property or not.

%
%
\section{A conjecture} \label{sec:conj}

Throughout this note, we assume $K$ is an arbitrary field unless otherwise specified.

We consider level Artinian monomial almost complete intersections in $R = K[x,y,z]$.  These are precisely the ideals of the form
\begin{equation} \label{e:levelmaci}
    I = (x^{\al + t}, y^{\be + t}, z^{\ga + t}, x^\al y^\be z^\ga).
\end{equation}
where $0 \leq \al \leq \be \leq \ga$ and $0 < t$, as shown in Section~6 of~\cite{MMN}.

Given known results and extensive computations, the authors of~\cite{MMN} made the following conjecture.
\begin{conjecture} \label{j:original}
    Let $K$ be an algebraically closed field of characteristic zero and let $I \subset R = K[x,y,z]$ be a level Artinian monomial almost complete intersection,
    that is, $I$ is as in~(\ref{e:levelmaci}).  Then:
    \begin{enumerate}
        \item $R/I$ has the weak Lefschetz property if any of the following conditions hold:
            \begin{enumerate}
                \item $\al = 0$,
                \item $\al + \be + \ga$ is not divisible by $3$,
                \item $\ga > 2(\al + \be)$, or
                \item $t < \frac{1}{3}(\al + \be + \ga)$.
            \end{enumerate}
        \item $R/I$ does not have the weak Lefschetz property if $(\al, \be, \ga, t)$ is $(2, 9, 13, 9)$ or $(3, 7, 14, 9)$.
        \item Assuming the parameters fail all conditions in~(i) and are not as in~(ii), then $R/I$ does not have the weak Lefschetz property if and only if
            $t$ is even and any of the following conditions hold:
            \begin{enumerate}
                \item $\al$ is even, $\al = \be$, and $\ga - \al \equiv 3 \ (\hbox{\rm mod } 6)$;
                \item $\al$ is odd, $\al = \be$, and $\ga - \al \equiv 0 \ (\hbox{\rm mod } 6)$; or
                \item $\al$ is odd, $\be = \ga$, and $\ga - \al \equiv 0 \ (\hbox{\rm mod } 3)$.
            \end{enumerate}
    \end{enumerate}
\end{conjecture}

Notice that the conditions in part~(iii) of Conjecture~\ref{j:original} can be restated in a more compact form.
\begin{conjecture} \label{j:restatement}
    Under the assumptions as in part~(iii) of Conjecture~\ref{j:original}, then $R/I$ does not have the weak Lefschetz property if and only if $t$ is even,
    $\al + \be + \ga$ is odd, and either $\al = \be$ or $\be = \ga$.
\end{conjecture}

In order to begin working on this conjecture the authors in~\cite{MMN} established a particular matrix in Theorem~7.2 and the corresponding Corollary~7.3,
whose determinant completely determines if the ring $R/I$ has the weak Lefschetz property.
\begin{theorem} \label{t:matrix-M}
    Let $K$ be an arbitrary field and let $I$ be as in~(\ref{e:levelmaci}) with the additional assumptions as in Conjecture~\ref{j:original}, part~(iii).
    Consider the square integer matrix $M$ of size $t + \frac{1}{3}(\al + \be - 2\ga)$:
    \[
        M =
        \left [ \begin{array}{cccccccccccccccccccccc}
            \binom{\ga}{\frac{1}{3}(\al+\be+\ga)} & \binom{\ga}{\frac{1}{3}(\al+\be+\ga) - 1} & \dots & \binom{\ga}{\ga-t+2} & \binom{\ga}{\ga-t+1} \\ \\
            \binom{\ga}{\frac{1}{3}(\al+\be+\ga)+1} & \binom{\ga}{\frac{1}{3}(\al+\be+\ga)} & \dots & \binom{\ga}{\ga-t+3} & \binom{\ga}{\ga-t+2} \\ \\
            && \vdots \\ \\
            \binom{\ga}{t-1} & \binom{\ga}{t-2} & \dots & \binom{\ga}{\frac{1}{3}(2\ga-\al-\be)+1}  & \binom{\ga}{\frac{1}{3}(2\ga-\al-\be)} \\ \\
            \binom{\ga+t}{t+\be-1} & \binom{\ga+t}{t+\be-2} & \dots & \binom{\ga+t}{\frac{1}{3}(2(\be+\ga)-\al)+1} & \binom{\ga+t}{\frac{1}{3}(2(\be+\ga)-\al)} \\ \\
            \binom{\ga+t}{t+\be-2} & \binom{\ga+t}{t+\be-3} & \dots & \binom{\ga+t}{\frac{1}{3}(2(\be+\ga)-\al)} & \binom{\ga+t}{\frac{1}{3}(2(\be+\ga)-\al)-1} \\ \\
            && \vdots \\ \\
            \binom{\ga+t}{t+\frac{1}{3}(\be+\ga-2\al)}& \binom{\ga+t}{t-1+\frac{1}{3}(\be+\ga-2\al)} & \dots & \binom{\ga+t}{\ga-\al+2} & \binom{\ga+t}{\ga-\al+1}
        \end{array} \right ].
    \]
    Then $\det{M} \equiv 0 \,( \hbox{\rm mod } \charf{K})$ if and only if $R/I$ fails to have the weak Lefschetz property.
\end{theorem}

Notice that the matrix $M$ has two distinct portions: a top half which has $t - \frac{1}{3}(\al+\be+\ga)$ rows and a bottom half which has
$\frac{1}{3}(2(\al+\be)-\ga)$ rows.  This will be especially useful in Section~\ref{sec:extremal-cases}.

A portion of Conjecture~\ref{j:original} has been proven; the results are summarised as follows:
\begin{remark} \label{r:conj-known}
    Part~(i) of Conjecture~\ref{j:original} is true by Corollary~6.3, Lemma~6.6, and Lemma~6.7 in~\cite{MMN}.  Part~(ii) is true by direct computation (e.g.,
    using a computer algebra system such as~\cite{cocoa} or~\cite{M2}).  Furthermore, the sufficiency of part~(iii) holds by Corollary~7.4 in~\cite{MMN}.
    Hence only the necessity of part~(iii) remains to be shown.
\end{remark}

Finally, we recall a result (\cite{MMN}, Proposition~6.1) that we make use of in Section~\ref{sec:unimodal}.
\begin{proposition} \label{p:resolution}
    Let $I$ be as in~(\ref{e:levelmaci}) and let $\sigma = \al + \be + \ga$.  Then $R/I$ has a free resolution of the form
    \begin{equation} \label{e:resolution}
        0
        \rightarrow
        R^3(-\sigma-2t)
        \rightarrow
        \begin{array}{c} R^3(-\sigma-t) \\ \oplus \\ R(-\al-\be-2t) \\ \oplus \\ R(-\al-\ga-2t) \\ \oplus \\ R(-\be-\ga-2t) \end{array}
        \rightarrow
        \begin{array}{c} R(-\sigma) \\ \oplus \\ R(-\al-t) \\ \oplus \\ R(-\be-t) \\ \oplus \\ R(-\ga-t) \end{array}
        \rightarrow
        R
        \rightarrow
        R/I
        \rightarrow
        0.
    \end{equation}
    Furthermore, if $\al > 0$ then this resolution is minimal.
\end{proposition}

%
%
\section{Some straightforward cases} \label{sec:simple-cases}

We establish necessary and sufficient numerical conditions for the weak Lefschetz property to hold in three families, all of which have the property $\al=\be=\ga$.

\begin{proposition} \label{p:abg-one}
    Suppose $\al=\be=\ga=1$ and $t \geq 1$.  Let $M$ be the matrix defined in Theorem~\ref{t:matrix-M}.  Then
    \[
        \det{M} = \left\{ \begin{array}{ll} 0, & \hbox{if $t$ is even;} \\ 2, & \hbox{if $t$ is odd.} \end{array} \right.
    \]
\end{proposition}
\begin{proof}
    Notice, $M \in \mathbb{Z}^{t \times t}$ is given by
    \[
        M = \left[
            \begin{array}{ccccccc}
                1 & 1 & 0 & \dots & 0 & 0 & 0 \\
                0 & 1 & 1 & \dots & 0 & 0 & 0 \\
                &&& \vdots &&& \\
                0 & 0 & 0 & \dots & 0 & 1 & 1 \\
                \binom{t+1}{t} & \binom{t+1}{t-1} & \binom{t+1}{t-2} & \dots & \binom{t+1}{3} & \binom{t+1}{2} & \binom{t+1}{1} \\
            \end{array}
        \right],
    \]
    which, after straightforward Gaussian elimination (\emph{Nota bene}: this requires only $t-1$ steps), yields the $t \times t$ matrix
    \[
        \tilde{M} = \left[
            \begin{array}{ccccccc}
                1 & 1 & 0 & \dots & 0 & 0 & 0 \\
                0 & 1 & 1 & \dots & 0 & 0 & 0 \\
                &&& \vdots &&& \\
                0 & 0 & 0 & \dots & 0 & 1 & 1 \\
                0 & 0 & 0 & \dots & 0 & 0 & \ell \\
            \end{array}
        \right]
    \]
    where
    \begin{eqnarray*}
        \ell    & = & \binom{t+1}{1} - \binom{t+1}{2} + \ldots + (-1)^{t+1}\binom{t+1}{t} \\
                & = & 1 - \sum_{i = 0}^{t+1}{(-1)^{i}\binom{t+1}{i}} + (-1)^{t+1} \\
                & = & 1 + (-1)^{t+1}. \\
    \end{eqnarray*}

    Hence, $\det{M} = \ell = 1 + (-1)^{t+1}$.
\end{proof}

Then the following is immediate using Theorem~\ref{t:matrix-M}.
\begin{corollary} \label{c:abg-one}
    Suppose $I = (x^{t+1}, y^{t+1}, z^{t+1}, xyz)$ where $t \geq 1$. Then the algebra $R/I$ has the weak Lefschetz property if and only if $t$ is odd and the
    characteristic of $K$ is not two.
\end{corollary}

\begin{proposition} \label{p:abg-two}
    Suppose $\al=\be=\ga=2$ and $t \geq 2$.  Let $M$ be the matrix defined in Theorem~\ref{t:matrix-M}.  Then
    \[
        \det{M} = \left\{ \begin{array}{ll} -t^2(t+3), & \hbox{if $t$ is even;} \\ (t+2)^2(t-1), & \hbox{if $t$ is odd.} \end{array} \right.
    \]
\end{proposition}
\begin{proof}
    Notice, $M \in \mathbb{Z}^{t \times t}$ is given by
    \[
        M = \left[
            \begin{array}{ccccccccc}
                1 & 2 & 1 & 0 & \dots & 0 & 0 & 0 & 0 \\
                0 & 1 & 2 & 1 & \dots & 0 & 0 & 0 & 0 \\
                &&&& \vdots &&&& \\
                0 & 0 & 0 & 0 & \dots & 0 & 1 & 2 & 1 \\
                \binom{t+2}{t+1} & \binom{t+2}{t} & \binom{t+2}{t-1} & \binom{t+2}{t-2} & \dots & \binom{t+2}{5} & \binom{t+2}{4} & \binom{t+2}{3} & \binom{t+2}{2} \\
                \binom{t+2}{t} & \binom{t+2}{t-1} & \binom{t+2}{t-2} & \binom{t+2}{t-3} & \dots & \binom{t+2}{4} & \binom{t+2}{3} & \binom{t+2}{2} & \binom{t+2}{1} \\
            \end{array}
        \right].
    \]
    We apply straightforward Gaussian elimination to all but the last row (\emph{Nota bene}: this requires only $2t-2$ steps) which yields the $t \times t$ matrix
    \[
        \tilde{M} = \left[
            \begin{array}{ccccccccc}
                1 & 2 & 1 & 0 & \dots & 0 & 0 & 0 & 0 \\
                0 & 1 & 2 & 1 & \dots & 0 & 0 & 0 & 0 \\
                &&&& \vdots &&&& \\
                0 & 0 & 0 & 0 & \dots & 0 & 1 & 2 & 1 \\
                0 & 0 & 0 & 0 & \dots & 0 & 0 & p & q \\
                0 & 0 & 0 & 0 & \dots & 0 & 0 & r & s \\
            \end{array}
        \right],
    \]
    where $p = (-1)^tt+t, q = (-1)^t(t-1) + 2t + 1, r = (-1)^t(t^2+t-1) + 1,$ and $s = (-1)^t(t^2 - 2) + 2.$  Then $\det{M} = \det{\tilde{M}} = p s - q r.$

    Suppose $t$ is even, then $$\det{M} = (2t)(t^2) - (3t)(t^2+t) = -t^2(t+3),$$ and if $t$ is odd, then $$\det{M} = (0)(-t^2+4) - (t+2)(-t^2-t+2) = (t+2)^2(t-1)$$
    as desired.
\end{proof}

\begin{corollary} \label{c:abg-two}
    Suppose $I = (x^{t+2}, y^{t+2}, z^{t+2}, x^2y^2z^2)$ where $t \geq 2$. Then the algebra $R/I$ has the weak Lefschetz property if the characteristic of $K$
    is zero or greater than $t+3$.
\end{corollary}

The following proposition is given without proof, as it directly imitates the proofs of Propositions~\ref{p:abg-one} and~\ref{p:abg-two}.
\begin{proposition} \label{p:abg-three}
    Suppose $\al=\be=\ga=3$ and $t \geq 3$.  Let $M$ be the matrix defined in Theorem~\ref{t:matrix-M}.  Then
    \[
        \det{M} = \left\{ \begin{array}{ll} 0, & \hbox{if $t$ is even;} \\ -\frac{1}{4}(t-1)^2(t+1)(t+2)(t+4)^2, & \hbox{if $t$ is odd.} \end{array} \right.
    \]
\end{proposition}

\begin{corollary} \label{c:abg-three}
    Suppose $I = (x^{t+3}, y^{t+3}, z^{t+3}, x^3y^3z^3)$ where $t \geq 3$. Then the algebra $R/I$ fails to have the weak Lefschetz property if $t$ is even.
    Further, $R/I$ has the weak Lefschetz property if $t$ is odd and either the characteristic of $K$ is zero or greater than $t+4$.
\end{corollary}

It is important to notice how the results in this section verify parts of Conjecture~\ref{j:restatement}:
\begin{remark} \label{r:abg-conj}
    For this remark, assume $K$ is a field of characteristic zero.

    In Corollaries~\ref{c:abg-one} and~\ref{c:abg-three} we have $\al + \be + \ga$ is odd, $\al = \be = \ga$, and $R/I$ has the weak Lefschetz
    property if and only if $t$ is odd.  This confirms Conjecture~\ref{j:restatement} for their respective cases.

    Further still, in Corollary~\ref{c:abg-two} we have that $\al + \be + \ga$ is even and $R/I$ always has the weak Lefschetz property.  This also confirms
    Conjecture~\ref{j:restatement} for the case $\al = \be = \ga = 2$.
\end{remark}

In the general case when $1 \leq \al = \be = \ga \leq t$, then the associated matrix $M$ defined in Theorem~$\ref{t:matrix-M}$ can be reduced to a matrix
$\tilde{M}$ of the form found in Propositions~\ref{p:abg-one} and~\ref{p:abg-two}.  That is, a diagonal matrix with entries $1$ on the diagonal except for
the bottom-right $\al \times \al$ matrix.  Hence, finding $\det{M}$ can be reduced to finding the determinant of an $\al \times \al$ matrix.

%
%
\section{Two extremal cases} \label{sec:extremal-cases}

In this section we consider two extremal cases for the parameters in Conjecture~\ref{j:restatement} where the weak Lefschetz property can be shown to hold.
We do this by computing the determinants of the associated matrices from Theorem~\ref{t:matrix-M}.

A nice concept that will allow a drastic simplification in the following determinants is the \emph{hyperfactorial}.
\begin{notation} \label{n:hyperfactorial}
    Let $n \geq 0$ be an integer.  Then define the \emph{hyperfactorial} of $n$ to be
    \[
        \HC(n) = \prod_{i=0}^{n-1}i!
    \]
    where it is important to notice that the product goes to $n-1$ and $\HC(0) = 1$.
\end{notation}

We need the following formula.
\begin{lemma} \label{l:binom-matrix-det}
    Let $T \geq B \geq 0$ be integers and let $N$ be an $n \times n$ matrix with entry $(i, j)$ given by
    \[
        N_{(i,j)} = \binom{T}{B+i-j} \quad \quad (1 \leq i, j \leq n).
    \]
    Then
    \[
        \det{N} = \frac{\HC(n)\HC(B)\HC(T-B)\HC(T+n)}{\HC(B+n)\HC(T-B+n)\HC(T)}.
    \]
\end{lemma}
\begin{proof}
    This follows by an application of Lemma~3 in~\cite{Kr} as described there on page~8.   We have written the result more conveniently, in particular, making
    use of the hyperfactorial form.
\end{proof}

We consider the case of $R/I$ as in Conjecture~\ref{j:restatement} where $\ga$ is maximal, that is, $\ga = 2(\al + \be)$.  Notice here, that the parameters
$\al, \be, \ga,$ satisfy the conditions of Theorem~\ref{t:matrix-M}.
\begin{theorem} \label{t:gamma-maximal}
    Let $1 \leq \al \leq \be, \ga = 2(\al + \be),$ and let $t \geq \frac{1}{3}(\al + \be + \ga) = \al + \be.$  Set $\de = t - (\al + \be)$.  Then the matrix
    $M$ from Theorem~\ref{t:matrix-M} is a $\de\times\de$ matrix which has entry $(i,j)$ given by
    \[
        M_{(i,j)} = \binom{\ga}{\al + \be + i - j} \quad \quad (1 \leq i, j \leq \de),
    \]
    and determinant
    \[
        \det{M} = \frac{\HC(\de)\HC^2(\al+\be)\HC(\ga + \de)}{\HC(\ga)\HC^2(t)}.
    \]
\end{theorem}
\begin{proof}
    First, notice that since $\ga = 2(\al + \be)$ the bottom half of $M$ from Theorem~\ref{t:matrix-M} has zero rows, so only the top half remains.  This
    gives precisely the matrix defined above.

    Hence, Lemma~\ref{l:binom-matrix-det} provides
    \begin{eqnarray*}
        \det{M} & = & \frac{\HC(\de)\HC(\al+\be)\HC(\ga - (\al + \be))\HC(\ga + \de)}{\HC(\al + \be + \de)\HC(\ga - (\al + \be) + \de)\HC(\ga)} \\
                & = & \frac{\HC(\de)\HC^2(\al+\be)\HC(\ga + \de)}{\HC(\ga)\HC^2(t)}
    \end{eqnarray*}
    where we use that $\ga - (\al + \be) = \al + \be$ and $\al + \be + \de = t$.
\end{proof}

As noted before, the parameters satisfy the conditions of Theorem~\ref{t:matrix-M}.
\begin{corollary} \label{c:gamma-maximal}
    Let $1 \leq \al \leq \be, \ga = 2(\al + \be),$ and let $t \geq \frac{1}{3}(\al + \be + \ga) = \al + \be.$  Consider the ideal given by
    \[
        I = (x^{\al+t}, y^{\be+t}, z^{\ga+t}, x^\al y^\be z^\ga) \subset R = K[x,y,z].
    \]
    Then $R/I$ has the weak Lefschetz property if $\charf{K} = 0$ or $\charf{K} \geq t + \al + \be$.
\end{corollary}
\begin{proof}
    Given the closed form of the determinant in Theorem~\ref{t:gamma-maximal}, it is clear that determinant is never zero.  Further still,
    we see that $\det{M}$ is not divisible by any prime equal to or greater than $t +  \al + \be = \ga + \delta$.
\end{proof}

We now consider the case of $R/I$ as in Conjecture~\ref{j:restatement} where $t$ is minimal, that is, $t = \frac{1}{3}(\al + \be + \ga).$  If we assume that
$1 \leq \al \leq \be \leq \ga \leq 2(\al + \be)$ and $\al + \be + \ga$ is divisible by $3$, then the parameters satisfy the conditions of
Theorem~\ref{t:matrix-M}.
\begin{theorem} \label{t:t-minimal}
    Let $1 \leq \al \leq \be \leq \ga \leq 2(\al+\be)$ such that $\al + \be + \ga$ is divisible by three and let $t = \frac{1}{3}(\al + \be + \ga).$  Set
    $\la = \frac{1}{3}(2(\al+\be)-\ga)$.  Then the matrix $M$ from Theorem~\ref{t:matrix-M} is a $\la\times\la$ matrix which has entry $(i,j)$ given by
    \[
        M_{(i,j)} = \binom{\ga+t}{\be+t+1-i-j} \quad \quad (1 \leq i, j \leq \la),
    \]
    and determinant
    \[
        \det{M} = (-1)^{\binom{\la}{2}} \frac{\HC(2t-\ga)\HC(2t-\be)\HC(2t-\al)\HC(\al+\be+\ga)}{\HC(\al+t)\HC(\be+t)\HC(\ga+t)}
    \]
\end{theorem}
\begin{proof}
    First, notice that since $t = \frac{1}{3}(\al+\be+\ga)$ then the top half of $M$ from Theorem~\ref{t:matrix-M} has zero rows, so only the bottom half
    remains.  This gives precisely the matrix defined above.

    In order to compute the determinant, we must first "flip" the matrix upside down.  This can be done in $\binom{\la}{2}$ operations (in can be done in less,
    but this gives a nice non-conditional form) yielding the matrix $\bar{M}$ such that $\det{M} = (-1)^{\binom{\la}{2}} \det{\bar{M}}$.  More
    importantly, the matrix $\bar{M}$ has for $1 \leq i,j \leq \la$ entry $(i,j)$ given by
    \[
        \bar{M}_{(i,j)} = \binom{\ga+t}{\be+t-\la+i-j}.
    \]
    Now we apply Lemma~\ref{l:binom-matrix-det} and obtain
    \begin{eqnarray*}
        \det{\bar{M}} & = & \frac{\HC(\la)\HC(\be + t -\la)\HC(\ga-\be+\la)\HC(\la + \ga + t)}{\HC(\be + t)\HC(\ga-\be+2\la)\HC(\ga + t)} \\
                      & = & \frac{\HC(2t-\ga)\HC(2t-\al)\HC(2t-\be)\HC(\al + \be + \ga)}{\HC(\be + t)\HC(\al+t)\HC(\ga+t)}
    \end{eqnarray*}
    where we use that $\ga + \la = 2t$ and $\al + \be = \la + t$.
\end{proof}

As noted before, the parameters satisfy the conditions of Theorem~\ref{t:matrix-M}.
\begin{corollary} \label{c:t-minimal}
    Let $1 \leq \al \leq \be \leq \ga \leq 2(\al+\be)$ such that $\al + \be + \ga$ is divisible by three and let $t = \frac{1}{3}(\al + \be + \ga).$  Consider
    the ideal given by
    \[
        I = (x^{\al+t}, y^{\be+t}, z^{\ga+t}, x^\al y^\be z^\ga) \subset R = K[x,y,z].
    \]
    Then $R/I$ has the weak Lefschetz property if $\charf{K} = 0$ or $\charf{K} \geq \al+\be+\ga$.
\end{corollary}
\begin{proof}
    Given the closed form of the determinant in Theorem~\ref{t:t-minimal}, it is clear that the determinant is never zero.  Further still, we see that $\det{M}$
    is not divisible by any prime equal to or greater than $\al + \be + \ga$ because $\al + \be + \ga-1$ is the maximum of the multiplicands in the numerator of
    the determinant.
\end{proof}

It is important to notice how the two results in this section verify parts of Conjecture~\ref{j:restatement}.
\begin{remark} \label{r:extremal-conj}
    Consider the case presented in Theorem~\ref{t:gamma-maximal}.  Notice that $\be < \ga$ and if $\al = \be$ then $\al + \be + \ga = 6\al$ is even, and this
    verifies Conjecture~\ref{j:restatement} in the case of $\ga$ being maximal.

    Consider now the case presented in Theorem~\ref{t:t-minimal}.  Notice that $t = \frac{1}{3}(\al + \be + \ga)$ is even if and only if $\al + \be + \ga$ is
    even.  Hence $t$ cannot be even at the same time as $\al + \be + \ga$ is odd, and this verifies Conjecture~\ref{j:restatement} for the case of $t$ being minimal.
\end{remark}

\begin{remark} \label{r:extremal-WLP}
    We notice that in the cases of $\ga$ being maximal and $t$ being minimal, the characteristics of $K$ where $R/I$ can possibly fail to have the weak Lefschetz
    property are bounded above by the maximum of the degrees of the generators of $I$.  However, in other cases described in Conjecture~\ref{j:restatement},
    this is not true.

    For example, consider the case $(\al, \be, \ga, t) = (2, 9, 13, 12)$ where the maximum degree of a generator of $I$ is $25$.  In this case,
    \begin{eqnarray*}
        \det{M} & = & -410893744849276115319750 \\
                & = & -2 \cdot 3^2 \cdot 5^3 \cdot 11^4 \cdot 13^5 \cdot 19 \cdot 23^3 \cdot 29 \cdot 5011.
    \end{eqnarray*}
    Hence, when $\charf{K} = 5011$ (or any other prime divisor of $\det{M}$) the algebra $R/I$ fails to have the weak Lefschetz property.
\end{remark}

There is a natural explanation why the determinants in Theorems~\ref{t:gamma-maximal} and~\ref{t:t-minimal} are non-trivial.  The determinants compute the
number of certain combinatorial objects.  More specifically, let $a, b, c$ be positive integers and consider a hexagon with side lengths $a, b, c, a, b, c$ with
angles $120^{\circ}$; a hexagon as described is called an \emph{$(a,b,c)$-hexagon}.  A \emph{lozenge} is a rhombus of unit side-length with angles $60^{\circ}$
and $120^{\circ}$.

\begin{figure}[h!] \label{f:hexagon}
    \includegraphics[scale=0.25]{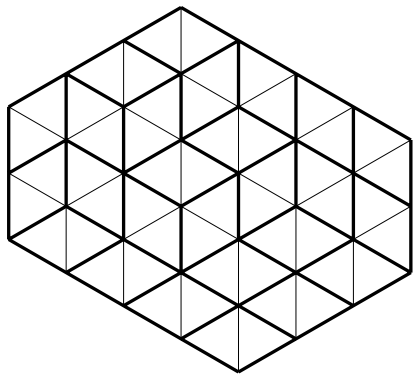}
    \caption{A $(2, 4, 3)$-hexagon tiled by lozenges.}
\end{figure}

The number of lozenge tilings is familiar (see Equation~(1.1) in ~\cite{CEKZ}).
\begin{proposition} \label{p:hexagon-tilings}
    Let $a, b, c \in \NN$.  Then the number of lozenge tilings
    of an $(a,b,c)$-hexagon is
    \[
        \frac{\HC(a)\HC(b)\HC(c)\HC(a+b+c)}{\HC(a+b)\HC(a+c)\HC(b+c)}.
    \]
\end{proposition}

Notice that if we set $a = n, b = B,$ and $c = T-B$ then the determinant found in Lemma~\ref{l:binom-matrix-det} counts the number of lozenge tilings of an
$(a,b,c)$-hexagon, i.e. an $(n, B, T-B)$-hexagon.  This connection is noted in both~\cite{CEKZ} and~\cite{Kr}.

In particular, the matrix $M$ associated to the case when $\ga$ is maximal (resp., $t$ is minimal) in Theorem~\ref{t:gamma-maximal} (resp.,
Theorem~\ref{t:t-minimal}) has determinant whose modulus counts the number of lozenge tilings of $(\al+\be, t-\al-\be, \al+\be)$-hexagons (resp.,
$(2t-\al, 2t-\be, 2t-\ga)$-hexagons).

This observation raises some natural questions:
\begin{enumerate}
    \item Can the connection above be extended in some way to all matrices $M$ appearing in Theorem~\ref{t:matrix-M}?
    \item More generally, can \emph{any} such combinatorial connection be found?
    \item Is there some natural property of the level Artinian monomial almost complete intersections which directly associates to the tilings of hexagons
          by lozenges?
\end{enumerate}

%
%
\section{The Hilbert function is peaked strictly unimodal} \label{sec:unimodal}

The Hilbert function is strongly tied to many properties of algebras.  If an algebra $A$ has the weak Lefschetz property, then its Hilbert function is
unimodal as shown in~\cite{HMNW}.  Ahn and Shin strengthened this result for level algebras.
\begin{proposition}[\cite{AS}, Theorem~3.6] \label{p:WLP-ps-unimodal}
    Let $A$ be a level Artinian standard graded $K$-algebra with the weak Lefschetz property.  Then the Hilbert function of $A$ is
    peaked strictly unimodal.
\end{proposition}

We have seen that some level Artinian monomial almost complete intersections in three variables fail to have the weak Lefschetz property.  Nevertheless, the
following result shows that their Hilbert functions are always peaked strictly unimodal regardless whether the quotient has the weak Lefschetz property or not.
\begin{lemma} \label{l:remainder-ps-unimodal}
    Let $1 \leq \al \leq \be \leq \ga < 2(\al + \be), t > \frac{1}{3}(\al + \be + \ga),$ and let $\al + \be + \ga$ be divisible by three.  Consider the ideal
    \[
        I = (x^{\al+t}, y^{\be+t}, z^{\ga+t}, x^\al y^\be z^\ga) \subset R = K[x,y,z].
    \]
    Then the Hilbert function of $R/I$ is peaked strictly unimodal with exactly two peaks in degrees $s = \frac{2}{3}(\al + \be + \ga) + t - 2$ and $s + 1$.
\end{lemma}
\begin{proof}
    Let $h$ be the Hilbert function of $R/I$, that is, $h(d) = \dim_K{[R/I]_d}$ for $d \in \ZZ$.

    By Proposition~\ref{p:resolution} the socle-degree of $R/I$ is $e = \al+\be+\ga+2t-3$ and the Cohen-Macaulay type of $R/I$ is three.  This implies
    $h(e) = 3$ and $h(d) = 0$ for $d > e$.

    Further, setting $s = \frac{2}{3}(\al + \be + \ga) + t - 2$, then, by Lemma~7.1 in~\cite{MMN}, $h(s) = h(s+1)$.  Moreover, since $2(\al + \be) > \ga$
    and as $\al + \be + \ga$ is divisible by three, then
    \begin{equation} \label{e:diff-abg}
        2(\al + \be) - \ga \geq 3.
    \end{equation}

    Now the proof is carried out in two steps.

    \emph{Step 1: Strict increase for $0 \leq d \leq s$.}
    First, notice that since $d \leq s$, then $\al + \be + \ga + t > s \geq d$ and further $\al + \be + 2t > \frac{1}{3}(4\al + 4\be + \ga) + t > s \geq d$.
    This implies that the ultimate and penultimate free modules in the Resolution~(\ref{e:resolution}) yield no contribution to the Hilbert function in degree $d$.  Hence,
    if $d \leq s$, we have that
    \begin{eqnarray*}
        h(d) & = & \binom{2+d}{2} - \binom{2+d-\al-\be-\ga}{2} \\
             &   & - \left[ \binom{2+d-\al-t}{2} + \binom{2+d-\be-t}{2} + \binom{2+d-\ga-t}{2} \right]
    \end{eqnarray*}
    and thus
    \begin{eqnarray*}
        h(d+1) - h(d) & = & \binom{2+d}{1} - \binom{2+d-\al-\be-\ga}{1} \\
                      &   & - \left[ \binom{2+d-\al-t}{1} + \binom{2+d-\be-t}{1} + \binom{2+d-\ga-t}{1} \right].
    \end{eqnarray*}

    For $0 \leq d < s$, when considering $h(d+1) - h(d)$, there are eight possible cases where the different binomial terms are nonzero in $h(d+1) - h(d)$.
    Furthermore, these eight cases are broken into two families:  when $d+1 < \al + \be + \ga$ and when $\al + \be + \ga \leq d+1$.

    Assume $d+1 < \al + \be + \ga$.
    \begin{enumerate}
        \item If $d + 1 < \al + t$, then
            \begin{eqnarray*}
                h(d+1) - h(d) & = & 2+d \\
                              & \geq & 2.
            \end{eqnarray*}
        \item If $\al + t \leq d + 1 < \be + t$, then
            \begin{eqnarray*}
                h(d+1) - h(d) & = & 2 + d - (2 + d - \al - t) \\
                              & = & \al + t \\
                              & \geq & 3.
            \end{eqnarray*}
        \item If $\be + t \leq d + 1 < \ga + t$, then
            \begin{eqnarray*}
                h(d+1) - h(d) & = & 2 + d - \left[ (2 + d - \al - t) + (2 + d - \be - t) \right] \\
                              & = & \al + \be + 2t - (2+d) \\
                              & \geq & \al + \be + t - \ga \\
                              & \geq & \frac{2}{3}(2(\al + \be) - \ga) + 1 \\
                              & \geq & 3,
            \end{eqnarray*}
            where the second inequality uses that $t > \frac{1}{3}(\al + \be + \ga)$, and the third
            inequality uses Inequality~(\ref{e:diff-abg}).
        \item If $\ga + t \leq d + 1$, then
            \begin{eqnarray*}
                h(d+1) - h(d) & = & 2 + d - \left[ (2 + d - \al - t) + (2 + d - \be - t) + (2 + d - \ga - t) \right] \\
                              & = & \al + \be + \ga + 3t - 2(2+d) \\
                              & \geq & \al + \be + \ga + 3t - 2(s+1) \\
                              & = & t - \frac{1}{3}(\al + \be + \ga) + 2 \\
                              & \geq & 3,
            \end{eqnarray*}
            where the second inequality uses that $t > \frac{1}{3}(\al + \be + \ga)$.
    \end{enumerate}

    Assume $\al + \be + \ga \leq d+1$.
    \begin{enumerate}
        \item If $d + 1 < \al + t$, then
            \begin{eqnarray*}
                h(d+1) - h(d) & = & 2+d - (2 + d - \al - \be - \ga) \\
                              & = & \al + \be + \ga \\
                              & \geq & 3.
            \end{eqnarray*}
        \item If $\al + t \leq d + 1 < \be + t$, then
            \begin{eqnarray*}
                h(d+1) - h(d) & = & 2+d - \left[ (2 + d - \al - \be - \ga) + (2 + d - \al - t) \right] \\
                              & = & 2\al + \be + \ga + t - (2+d) \\
                              & \geq & 2\al + \ga \\
                              & \geq & 3.
            \end{eqnarray*}
        \item If $\be + t \leq d + 1 < \ga + t$, then
            \begin{eqnarray*}
                h(d+1) - h(d) & = & 2+d - \left[ (2 + d - \al - \be - \ga) + (2 + d - \al - t) + (2 + d - \be - t) \right] \\
                              & = & 2\al + 2\be + \ga + 2t - 2(2+d) \\
                              & \geq & 2\al+2\be - \ga \\
                              & \geq & 3,
            \end{eqnarray*}
            where the second inequality uses Inequality~(\ref{e:diff-abg}).
        \item If $\ga + t \leq d + 1$, then
            \begin{eqnarray*}
                h(d+1) - h(d) & = & 2+d - \left[ (2 + d - \al - \be - \ga) + (2 + d - \al - t) \right] \\
                              &   & - \left[ (2 + d - \be - t) + (2 + d - \ga - t) \right] \\
                              & = & 2(\al + \be + \ga) + 3t - 3(2+d) \\
                              & \geq & 2(\al + \be + \ga) + 3t - 3(s+2) + 3 \\
                              & = & 3.
            \end{eqnarray*}
    \end{enumerate}

    Thus we have that $h(d+1) - h(d) > 0$ for all $0 \leq d < s$ implying that the Hilbert function is strictly increasing from degree $0$ to degree $s$.

    \emph{Step 2: Strict decrease for $s + 1\leq d \leq e.$}
    Let $k$ be the Hilbert function of the $K$-dual of $R/I$, that is of $(R/I)^{\vee}$.  Then $k(d) = \dim_K{[(R/I)^{\vee}]_d}$, so
    $h(d) = k(-d)$ for all $d \in \ZZ$.

    Since $d \geq s+1$, then $\al + \be + \ga \leq s + 1 \leq d$ and, using Inequality~(\ref{e:diff-abg}), $\al + t \leq \be + t \leq \ga + t \leq s + 1 \leq d$.
    This implies that the ultimate and penultimate free
    modules in the resolution of $(R/I)^{\vee}$ (which is dual to the Resolution~(\ref{e:resolution})) yield no contribution to the
    Hilbert function of $(R/I)^{\vee}$ in degree $-d$.  Hence, if $s+1 \leq d$, we have that
    {\footnotesize
        \begin{eqnarray*}
            k(-d) & = & 3\binom{-d-1+\al+\be+\ga+2t}{2} - 3\binom{-d-1+\al+\be+\ga+t}{2} \\
                 &   & - \left[ \binom{-d-1+\be+\ga+2t}{2} + \binom{-d-1+\al+\ga+2t}{2} + \binom{-d-1+\al+\be+2t}{2} \right]
        \end{eqnarray*}
    }
    and thus
    {\footnotesize
        \begin{eqnarray*}
            k(-d) - k(-d-1) & = & 3\binom{-d-2+\al+\be+\ga+2t}{1} - 3\binom{-d-2+\al+\be+\ga+t}{1} \\
                        &   & - \left[ \binom{-d-2+\be+\ga+2t}{1} + \binom{-d-2+\al+\ga+2t}{1} + \binom{-d-2+\al+\be+2t}{1} \right].
        \end{eqnarray*}
    }

    For $s+1 \leq d < e$, when considering $h(d) - h(d+1) = k(-d) - k(-d-1)$, there are eight possible cases where the different binomial terms are nonzero in
    $k(-d)-k(-d-1)$.  Furthermore, these eight cases are broken into two families:  when $d + 1 \leq \al+\be+\ga+t-2$ and when $\al+\be+\ga+t-2 < d+1$.

    Assume $\al + \be + \ga + t-2 < d+1$.
    \begin{enumerate}
        \item If $\be + \ga + 2t - 2 < d+1$, then
            \begin{eqnarray*}
                k(-d) - k(-d-1) & = & 3(-d-2+\al+\be+\ga+2t) \\
                                & \geq & 3(-(e-1)-2+\al+\be+\ga+2t) \\
                                & = & 6,
            \end{eqnarray*}
            where we use $d+1 \leq e = \al + \be + \ga + 2t - 3$.
        \item If $\al + \ga + 2t - 2 < d+1 \leq \be + \ga + 2t - 2$, then
            \begin{eqnarray*}
                k(-d) - k(-d-1) & = & 3(-d-2+\al+\be+\ga+2t) - (-d-2+\be+\ga+2t) \\
                                & = & 3\al + 2\be + 2\ga + 4t - 4 - 2d \\
                                & \geq & 3\al + 2 \\
                                & \geq & 5.
            \end{eqnarray*}
        \item If $\al + \be + 2t - 2 < d+1 \leq \al + \ga + 2t - 2$, then
            \begin{eqnarray*}
                k(-d) - k(-d-1) & = & 3(-d-2+\al+\be+\ga+2t) - (-d-2+\be+\ga+2t) \\
                                &   & - (-d-2+\al+\ga+2t) \\
                                & = & 2\al + 2\be + \ga + 2t - 2 - d \\
                                & \geq & \al + 2\be + 1 \\
                                & \geq & 4.
            \end{eqnarray*}
        \item If $s + 1 < d+1 \leq \al + \be + 2t - 2$, then
            \begin{eqnarray*}
                k(-d) - k(-d-1) & = & 3(-d-2+\al+\be+\ga+2t) - (-d-2+\be+\ga+2t) \\
                                &   & - (-d-2+\al+\ga+2t) - (-d-2+\al+\be+2t) \\
                                & = & \al + \be + \ga \\
                                & \geq & 3.
            \end{eqnarray*}
    \end{enumerate}

    Assume $\al + \be + \ga + t-2 < d+1$.
    \begin{enumerate}
        \item If $\be + \ga + 2t - 2 < d+1$, then
            \begin{eqnarray*}
                k(-d) - k(-d-1) & = & 3(-d-2+\al+\be+\ga+2t) - 3(-d-2+\al+\be+\ga+t) \\
                                & = & 3t \\
                                & \geq & 6.
            \end{eqnarray*}
        \item If $\al + \ga + 2t - 2 < d+1 \leq \be + \ga + 2t - 2$, then
            \begin{eqnarray*}
                k(-d) - k(-d-1) & = & 3(-d-2+\al+\be+\ga+2t) - 3(-d-2+\al+\be+\ga+t) \\
                                &   & - (-d-2+\be+\ga+2t) \\
                                & = & t + d + 2 - \be - \ga \\
                                & \geq & 3t + \al - \be \\
                                & \geq & 3 + 2\al + \ga \\
                                & \geq & 6,
            \end{eqnarray*}
            where the second inequality uses that $t > \frac{1}{3}(\al+\be+\ga)$.
        \item If $\al + \be + 2t - 2 < d+1 \leq \al + \ga + 2t - 2$, then
            \begin{eqnarray*}
                k(-d) - k(-d-1) & = & 3(-d-2+\al+\be+\ga+2t) - 3(-d-2+\al+\be+\ga+t) \\
                                &   & - (-d-2+\be+\ga+2t) - (-d-2+\al+\ga+2t) \\
                                & = & 2d + 4 - \al - \be - 2\ga - t \\
                                & \geq & \al + \be - 2\ga + 3t \\
                                & \geq & 2\al + 2\be - \ga + 3 \\
                                & \geq & 6,
            \end{eqnarray*}
            where the second inequality uses that $t > \frac{1}{3}(\al+\be+\ga)$, and the third inequality uses Inequality~(\ref{e:diff-abg}).
        \item If $s + 1 < d+1 \leq \al + \be + 2t - 2$, then
            \begin{eqnarray*}
                k(-d) - k(-d-1) & = & 3(-d-2+\al+\be+\ga+2t) - 3(-d-2+\al+\be+\ga+t) \\
                                &   & - (-d-2+\be+\ga+2t) - (-d-2+\al+\ga+2t) \\
                                &   & - (-d-2+\al+\be+2t) \\
                                & = & 3d + 6 - 2(\al+\be+\ga)-3t \\
                                & \geq & 3.
            \end{eqnarray*}
    \end{enumerate}

    Hence we have that $h(d) - h(d+1) = k(-d) - k(-d-1) > 0$ for all $s+1 \leq d < e$ implying that the Hilbert function is decreasing from $s+1$ to the
    socle degree $e$.
\end{proof}

This provides the following result which gives an affirmative answer to (part of) Question~8.2(1) in~\cite{MMN}.
\begin{theorem} \label{t:maci-ps-unimodal}
    Let $I \subset R = K[x,y,z]$ be a level Artinian monomial almost complete intersection.   Then $R/I$ has a peaked strictly unimodal Hilbert function.
\end{theorem}
\begin{proof}
    In case~(i) of Conjecture~\ref{j:original}, Remark~\ref{r:conj-known} guarantees the weak Lefschetz property of $R/I$, so the claim follows by
    Proposition~\ref{p:WLP-ps-unimodal}.

    Similarly, we conclude in case~(iii) of Conjecture~\ref{j:original}, when $\ga$ is maximal or $t$ is minimal by using Corollaries~\ref{c:gamma-maximal}
    and~\ref{c:t-minimal}.

    In all the remaining cases of Conjecture~\ref{j:original}, we conclude by Lemma~\ref{l:remainder-ps-unimodal}.
\end{proof}

%
%

\end{document}